\documentclass[12pt,a4paper]{article}
\usepackage[latin1]{inputenc}
\usepackage[french]{babel}
\usepackage[T1]{fontenc}
\usepackage{amssymb,theorem}
\usepackage{graphics}
\newtheorem{thm}{Th\'eor\`eme}

\newtheorem{prop}[thm]{Proposition}

 {\theorembodyfont{\rmfamily}
\newtheorem{exemp}[thm]{Exemple}}
\newcommand{\R}{\mathbb R}
\newcommand{\C}{\mathbb C}

\newcommand{\N}{\mathbb N}
\newcommand{\Z}{\mathbb Z}
\newcommand{\Q}{\mathbb Q}

\begin{document}
\title{\Large{\textbf{Relations de r\'ecurrence linéaires, primitivit\'e et loi de Benford}}}
\author{Hugues Deligny et Paul Jolissaint}
\maketitle

\section{Introduction}

Nous nous int\'eressons ici \`a la loi de Benford pour des suites $(a_n)_{n\geq 1}\subset\R^+ $ qui satisfont une relation de r\'ecurrence d'ordre $k$ de la forme :
$$
(\star)\quad a_{n+k}=c_{k-1}a_{n+k-1}+c_{k-2}a_{n+k-2}+\ldots+c_1a_{n+1}+c_0a_n
$$
avec les conditions initiales $a_i>0$ pour tout $1\leq i\leq k$. On supposera la plupart du temps que les coefficients $c_i$ sont positifs ou nuls, et que $c_0>0$, de sorte que la relation de r\'ecurrence soit effectivement d'ordre $k$. On va pr\'esenter des conditions suffisantes sur les $c_i$ pour qu'il existe (au moins) un entier $b>2$ tel que la suite $(a_n)_{n\geq 1}$ satisfasse la loi de Benford en base $b$.
\par
Rappelons qu'\'etant donn\'e un tel entier $b>2$, une suite $(u_n)_{n\geq 1}\subset\R^+\setminus\{0\}$ satisfait la \emph{loi de Benford en base} $b$ si, pour tout $t\in [1,b)$, on a :
$$
\lim_{N\to\infty}\frac{|\{1\leq n\leq N\ :\ M_b(u_n)<
  t\}|}{N}=\log_b(t),
$$
où $M_b(u_n)$ d\'esigne la \emph{mantisse} de $u_n$, c'est-\`a-dire l'unique \'el\'ement de $[1,b)$ tel que $u_n=M_b(u_n)\cdot b^m$ avec $m\in\Z$
(cf \cite{jol2}, d\'efinition 2.1). La d\'efinition ci-dessus g\'en\'eralise le cas classique qui correspond \`a $b=10$ (dix) et qui affirme que le premier chiffre significatif $d_1(u_n)$ de $u_n$ satisfait la loi de probabilit\'e (cf remarque (2) de \cite{jol2}) :\\ 
pour tout chiffre $1\leq d\leq 9$, on a 
$$
\lim_{N\to\infty}\frac{|\{1\leq n\leq N:d_1(u_n)=d\}|}{N}=\log_{10}\left(1+\frac{1}{d}\right).
$$
\par\vspace{3mm}
D'apr\`es les th\'eor\`emes principaux de \cite{jol} et \cite{jol2}, une suite $(u_n)\subset\R^+\setminus\{0\}$ satisfait la loi de Benford en base $b$ d\`es qu'elle remplit les deux conditions suivantes :
\begin{enumerate}
	\item [(a)] il existe des nombres r\'eels $\alpha>0$, $\rho>0$ et $\mu$ tels que $\displaystyle{\lim_{n\to\infty} \frac{u_n}{n^\mu\rho^n}=\alpha}$ ;
	\item [(b)] $\log_b(\rho)$ est irrationnel.
\end{enumerate}
\par
Si c'est le cas, si $Q(x)$ est un polynôme non constant \`a coefficients entiers et $Q(x)\geq 1$ pour $x$ assez grand,
la sous-suite $(u_{Q(n)})_{n\geq 1}$ satisfait \'egalement la loi de Benford en base $b$.
\par\vspace{3mm}
L'exemple suivant découle imm\'ediatement des deux conditions ci-dessus, et il ne semble pas avoir été découvert auparavant :


\begin{exemp} Soit $(p_n)_{n\geq 1}$ la suite croissante des nombres premiers. Comme cela a \'et\'e observ\'e dans \cite{dia}, la suite $(p_n)$ ne satisfait pas la loi de Benford. En revanche, choisissons deux entiers $\ell\geq 2$ et $b>2$ tels que $\log_b(\ell)$ soit irrationnel. Alors la suite $(p_{\ell^n})$ satisfait la loi de Benford en base $b$. En effet, par le th\'eor\`eme des nombres premiers, on a 
$$
\lim_{n\to\infty}\frac{p_n}{n\ln(n)}=1.
$$
En passant \`a la sous-suite $(\ell^n)_{n\geq 1}\subset\N^*$, on a 
$$
\lim_{n\to\infty}\frac{p_{\ell^n}}{n\ell^n\ln(\ell)}=1\quad\textrm{et\ donc}\quad
\lim_{n\to\infty}\frac{p_{\ell^n}}{n\ell^n}=\ln(\ell).
$$
\end{exemp}
\par\vspace{3mm}
Dans le cas d'une suite $(a_n)_{n\geq 1}$ qui satisfait une relation de r\'ecurrence du type $(\star)$, $(a_n)_{n\geq 1}$ s'exprime \`a l'aide des racines du \emph{polyn\^ome caract\'eristique} de la relation de r\'ecurrence, c'est-\`a-dire le polyn\^ome 
$$
p(x)=x^k-c_{k-1}x^{k-1}-\ldots-c_1x-c_0.
$$
Plus pr\'ecis\'ement, \'ecrivons $p(x)=(x-\zeta_1)^{\mu_1}\cdots (x-\zeta_m)^{\mu_m}$ où $\zeta_1,\ldots,\zeta_m\in\C$ sont les z\'eros distincts de $p(x)$ de multiplicit\'es respectives $\mu_1,\ldots,\mu_m\geq 1$. Alors $(a_n)_{n\geq 1}$ est une combinaison lin\'eaire des $k$ suites $(n^\ell \zeta_j^n)_{n\geq 1}$ pour $0\leq\ell<\mu_j$ et $1\leq j\leq m$ : il existe des constantes $\alpha_{j,\ell}$ d\'ependant des conditions initiales $a_1,\ldots,a_k$ telles que
$$
a_n=\sum_{j=1}^m\sum_{\ell=0}^{\mu_j-1}\alpha_{j,\ell}n^\ell \zeta_j^n\quad \forall n\geq 1.
$$
D\`es lors, si $p(x)$ admet une racine positive simple $\rho$ telle que $|\zeta|<\rho$ pour toute autre racine $\zeta$ de $p(x)$ et si le coefficient de $\rho^n$ est positif dans la d\'ecomposition ci-dessus, alors la suite $(a_n)_{n\geq 1}$ satisfait la condition (a).
\par\vspace{3mm}
Il nous a sembl\'e int\'eressant de donner des conditions suffisantes sur les \emph{coefficients} de la relation de r\'ecurrence $(\star)$ plut\^ot que sur les \emph{racines} du polyn\^ome caract\'eristique de la relation de r\'ecurrence. Nous mettrons plus particulièrement l'accent sur le cas où les coefficients $c_i$ sont positifs ou nuls.

\par\vspace{3mm}
Le premier r\'esultat de notre \'etude est :

\begin{thm}
Soient $c_0,c_1,\ldots,c_{k-1}$ des nombres r\'eels. On pose 
$$
I=\{1\leq j\leq k-1:c_j\not=0\},
$$
et on suppose que $I\not=\emptyset$. Enfin, on associe à la suite $c_0,\ldots,c_{k-1}$ le polynôme 
$$
p(x)=x^k-c_{k-1}x^{k-1}-\ldots -c_1x-c_0.
$$ 
\begin{enumerate}
	\item [(1)] Si les $c_i$ sont positifs ou nuls et si $c_0>0$, le polyn\^ome $p$
admet une unique racine positive $\rho$ et $|\zeta|<\rho$ pour tout autre zéro $\zeta$ de $p$ (compt\'e avec multiplicit\'e) si, et seulement si $\mathrm{pgcd}(I\cup\{k\})=1$.
\item [(2)] Si $c_{k-1}\geq 2$ et si $c_{k-1}>\sum_{j=0}^{k-2}|c_j|+1$, alors le polyn\^ome $p$ 
admet une racine simple $c_{k-1}-1<\rho<c_{k-1}+1$, et toute autre racine $\zeta$ de $p$ satisfait $|\zeta|<1$.
\item [(3)] Supposons que les $c_i$ satisfassent l'une des deux conditions ci-dessus,
soit $b>2$ un entier tel que $\log_b(\rho)$ soit irrationnel, et soit une suite $(a_n)_{n\geq 1}\subset\R^+\setminus\{0\}$ qui satisfait la relation de r\'ecurrence 
$$
a_{n+k}=c_{k-1}a_{n+k-1}+\ldots+c_1a_{n+1}+c_0a_n.
$$
Si les $c_i$ satisfont la condition (1) ou la condition (2), et dans ce dernier cas si de plus $(a_n)$ ne tend pas vers 0, alors
la suite $(a_n)_{n\geq 1}$ satisfait la loi de Benford en base $b$.
Enfin, il en est de même de toute sous-suite de la forme $(a_{Q(n)})$ où $Q(x)$ d\'esigne un polyn\^ome non constant \`a coefficients entiers tel que $Q(x)\geq 1$ pour tout $x$ suffisamment grand.
\end{enumerate}
\end{thm}

\par\vspace{3mm}\noindent
\textbf{Remarque.} La condition supplémentaire $a_n\not\to 0$ dans la troisième affirmation du théorème ci-dessus est indispensable. En effet, la suite $a_n=10^{-n}$ tend vers 0, satisfait la relation de récurrence
$$
a_{n+2}=\frac{31}{10}a_{n+1}-\frac{3}{10}a_n
$$
dont le polynôme caractéristique est $p(x)=(x-3)(x-1/10)$ ; ses
coefficients satisfont la condition (2) du théorème, mais la suite ne satisfait pas la loi de Benford pour $b=10$, bien que $\log_{10}(3)$ soit irrationnel.

\par\vspace{3mm}
On verra que lorsque les $c_i$ sont positifs ou nuls, le polyn\^ome caract\'eristique $p(x)=x^k-c_{k-1}x^{k-1}-\ldots -c_1x-c_0$ admet toujours une racine dominante $\rho>0$, c'est-\`a-dire telle que $|\zeta|\leq\rho$ pour toute autre racine $\zeta$ ; regardons maintenant le cas où il admet
$h>1$ racines $\zeta_1=\rho,\ldots,\zeta_h$ telles que $|\zeta_j|=\rho$.
On note encore
$I=\{1\leq j<k:c_j>0\}$, et pour tout $m\in\{0,\ldots,h-1\}$, on pose 
$$
I_m=\{j\in I : j\equiv m (\bmod\ h)\}
$$
de sorte que $(I_m)_{0\leq m<h}$ constitue une partition de $I$.

\begin{thm}
Soient $c_0>0, c_1\geq 0,\ldots,c_{k-1}\geq 0$, et $h$ comme ci-dessus. Alors $h=\mathrm{pgcd}(I\cup\{k\})$. De plus, toute suite $(a_n)_{n\geq 1}$ qui satisfait la relation de r\'ecurrence
$$
a_{n+k}=c_{k-1}a_{n+k-1}+c_{k-2}a_{n+k-2}+\ldots+c_1a_{n+1}+c_0a_n
$$
avec des conditions initiales positives ou nulles $a_1,\ldots,a_k$ est r\'eunion des $h$ sous-suites $((a_{m+hn})_{n\geq 1})_{0\leq m<h}$. Enfin, pour tout $0\leq m<h$, on a l'alternative suivante : 
\begin{itemize}
	\item Si $\{j\in I_m:a_j>0\}\not=\emptyset$, alors $\lim_{n\to\infty}\frac{a_{m+hn}}{\rho^n}$ existe et est positive, et la suite $(a_{m+hn})_{n\geq 1}$ satisfait la loi de Benford en base $b$ pour tout $b$ tel que $\log_b(\rho)\notin\Q$.
	\item Si $\{j\in I_m:a_j>0\}=\emptyset$, alors $a_{m+hn}=0$ pour tout $n\geq 1$.
\end{itemize}
\end{thm}


\begin{exemp} Le th\'eor\`eme 2 s'applique au cas des \emph{suites de Fibonacci d'ordre} $k\geq 2$ : on choisit des nombres r\'eels $a_1,\ldots,a_k>0$ arbitraires, et on d\'efinit $(a_n)_{n\geq 1}$ par
$$
a_{n+k}=a_{n+k-1}+a_{n+k-2}+\ldots+a_{n+1}+a_n.
$$
\end{exemp}
Notons que la racine dominante $\rho>0$ est irrationnelle grâce \`a l'observation suivante : 

\begin{prop}
Soit $k\geq 2$ un entier et soient $c_1,\ldots,c_{k-1}\in\N$, et $p(x):=x^k-c_{k-1}x-\ldots-c_1x-1$. Si $\rho>0$ est une racine rationnelle de $p(x)$, alors $\rho=1$.
\end{prop}
\emph{Preuve.} \'Ecrivons $\rho=\frac{p}{q}$ avec $p,q\in\N^*$ et $\mathrm{pgcd}(p,q)=1$. Alors, en utilisant le fait que $\rho$ est une racine de $p(x)$, on obtient
$$
p^k=c_{k-1}p^{k-1}q+c_{k-2}p^{k-2}q^2+\ldots+c_1pq^{k-1}+q^k
$$
qui implique que $p$ divise $q$.
\hfill $\square$

\par\vspace{3mm}
Enfin, il est n\'ecessaire de pouvoir pr\'eciser les valeurs de $b$ pour lesquelles la loi de Benford est satisfaite au moins lorsque les coefficients $c_j$ sont rationnels. 

\begin{thm}
Soit $(a_n)_{n\geq 1}\subset\R^+\setminus\{0\}$ une suite qui satisfait la relation de r\'ecurrence d'ordre $k$ 
$$
a_{n+k}=c_{k-1}a_{n+k-1}+c_{k-2}a_{n+k-2}+\ldots+c_1a_{n+1}+c_0a_n
$$ 
avec des coefficients $c_i\geq 0$, $c_0>0$ et $\{i<k:c_i>0\}\not=\emptyset$, et avec
les valeurs initiales $a_1,\ldots,a_k$ strictement positives. Alors :
\begin{enumerate}
	\item [(1)] La suite $(a_n)_{n\geq 1}$ satisfait la loi de Benford dans presque toute base $b$ au sens suivant :
	$$
	\lim_{N\to\infty}\frac{1}{N}|\{2<b\leq N:(a_n)\ \textrm{ne\ suit\ pas\ la\ loi\ de\ Benford\ en\ base}\ b\}|=0.
	$$
	\item [(2)] Si les coefficients $c_i$ sont tous rationnels et si la racine r\'eelle positive du polyn\^ome caract\'eristique de la relation de r\'ecurrence n'est ni un entier ni un inverse d'entier, alors $(a_n)_{n\geq 1}$ satisfait la loi de Benford dans toute base $b>2$.
	\end{enumerate}
\end{thm}

\par\vspace{3mm}\noindent
\textbf{Remarque.} Dans le cas o\`u seul $c_0>0$, c'est-à-dire si $I=\emptyset$, toute suite $(a_n)_{n\geq 1}$ qui satisfait la relation de récurrence $(\star)$ est de la forme :
$$
a_n=(\sqrt[k]{c_0})^n\cdot\sum_{\ell=0}^{k-1}\alpha_\ell e^{2\pi i\ell n/k}\quad\forall n.
$$
Par p\'ediodicit\'e de $e^{2\pi i\ell n/k}$, la suite $(a_n)$ est donc r\'eunion de $k$ sous-suites $(a_{m+kn})_{n\geq 1}$, $0\leq m<k$, qui sont toutes g\'eom\'etriques de raison $c_0$.
\par
Ainsi, si $c_0\not=1$, chaque sous-suite $(a_{m+kn})_{n\geq 1}$ satisfait la loi de Benford dans presque toute base. Si en revanche $c_0=1$, chaque sous-suite est constante et ne satisfait pas la loi de Benford.

\par\vspace{3mm}
Les preuves utilisent principalement la matrice compagnon du polyn\^ome $p(x)$ et elles
reposent sur la th\'eorie de Perron-Frobenius des matrices \`a coefficients positifs ou nuls qui sont irr\'eductibles ou primitives, et dont nous rappelons les principaux r\'esultats dans le paragraphe suivant. Les preuves des th\'eor\`emes 2 et 3 se trouvent dans le paragraphe 3 et la preuve du th\'eor\`eme 6 dans le dernier paragraphe.

\section{Matrices irr\'eductibles ; matrices primitives}

Nous rappelons ci-dessous les d\'efinitions et les r\'esultats principaux de la th\'eorie de Perron-Frobenius \`a propos des matrices irr\'eductibles et primitives. Nos r\'ef\'erences sont d'une part le chapitre 8 de la monographie de Carl D. Meyer \cite{meyer} et d'autre part les notes de J E Rombaldi \cite{romb}.
\par\vspace{3mm}
Soit $C\in M_k(\R)$ une matrice $k\times k$ \`a coefficients positifs ou nuls (on note $C\geq 0$). Nous consid\'ererons ici les valeurs propres complexes de $C$, c'est-\`a-dire ses valeurs propres en tant qu'endomorphisme de $\C^k$. L'ensemble des valeurs propres est le \emph{spectre} de $C$ et sera not\'e $\sigma(C)$.
Nous rappelons que le \emph{rayon spectral} de $C$ est 
$$
\rho(C)=\max\{|\lambda|:\lambda\in\sigma(C)\}.
$$
On dit que $C$ est \emph{r\'eductible} s'il existe une matrice de permutation $P$ telle que $PCP^{-1}$ soit de la forme :
$$
PCP^{-1}=
\left(
\begin{array}{cc}
C' & C'''\\
0 & C''
\end{array}\right)
$$
où $C'$ et $C''$ sont des matrices carr\'ees de dimensions positives, et on dit qu'elle est \emph{irr\'eductible} si elle n'est pas r\'eductible.
\par
On d\'emontre qu'une telle matrice $C=(c_{i,j})$ est irr\'eductible si et seulement si le \emph{graphe orient\'e} $G(C)$ associ\'e est fortement connexe. Le graphe $G(C)$ a pour ensemble de sommets les entiers $1,2\ldots,k$, et il y a une arête \emph{orient\'ee} de $i$ vers $j$ si et seulement si $c_{i,j}>0$ ; le graphe est \emph{fortement connexe} s'il existe un chemin orient\'e de 1 vers 1 passant par tous les sommets. 
\par
Si $C$ est irr\'eductible, le rayon spectral $\rho=\rho(C)$ est strictement positif et c'est une valeur propre simple du polyn\^ome caract\'eristique $p_C(x)=\det(C-xI)$. Le sous-espace propre correspondant est de dimension (complexe) \'egale \`a 1, et il est engendr\'e par un \emph{vecteur de Perron-Frobenius}, c'est-\`a-dire un vecteur 
$v=\left(
\begin{array}{c}
v_1\\
\vdots\\
v_k\end{array}\right)$
dont toutes les composantes $v_i$ sont strictement positives.
\par
Notons $h\geq 1$ le nombre de valeurs propres $\lambda\in\sigma(C)$ telles que $|\lambda|=\rho$. Si $h>1$, on l'appelle \emph{l'indice d'imprimitivit\'e} de $C$, et 
si $h=1$, on dit que $C$ est \emph{primitive}.
\par
On d\'emontre que si $C$ est irr\'eductible, alors elle est primitive si et seulement s'il existe un entier $m>0$ tel que $C^m>0$, c'est-\`a-dire les coefficients de $C^m$ sont tous strictement positifs (test de primitivit\'e de Frobenius, cf \cite{meyer}, pp. 674 et 678).
\par
Si $C$ est irr\'eductible mais non primitive, si $\{\lambda_1,\ldots,\lambda_h\}$ est l'ensemble des valeurs propres de $C$ de module \'egal \`a $\rho$, alors 
$$
\{\lambda_1,\ldots,\lambda_h\}=\{\rho,\rho\omega,\rho\omega^2,\ldots,\rho\omega^{h-1}\}
$$
où $\omega=e^{2\pi i/h}$ (\cite{meyer}, p. 676).
\par
Soit $S=(s_{i,j})\in M_k(\R)$ une matrice \`a coefficients positifs ou nuls. On dit que $S$ est \emph{stochastique} (par rapport aux lignes) si, pour tout $1\leq i\leq k$, on a
$$
\sum_{j=1}^k s_{i,j}=1.
$$
Si une telle matrice $S$ est de plus irr\'eductible, si $h$ est son indice d'imprimitivit\'e, alors elle admet toutes les racines $h$-i\`eme de l'unit\'e comme valeurs propres.
\par
Grâce au test de primitivit\'e de Frobenius, on observe que si $A$ et $B$ sont des matrices \`a coefficients positifs ou nuls irr\'eductibles et s'il existe $\epsilon>0$ tel que $A\geq \epsilon B$ (c'est-\`a-dire si $a_{i,j}\geq\epsilon b_{i,j}$ pour tous $i,j$), alors $A$ est primitive si $B$ l'est.
\par
Consid\'erons encore une matrice $C\geq 0$ irr\'eductible et notons 
$v=\left(
\begin{array}{c}
v_1\\
\vdots\\
v_k\end{array}\right)>0$ un vecteur de Perron-Frobenius de $C$ ;
la transpos\'ee $C^T$ de $C$ poss\`ede les mêmes propri\'et\'es, donc elle admet un vecteur de Perron-Frobenius 
$w=\left(
\begin{array}{c}
w_1\\
\vdots\\
w_k\end{array}\right)>0$ associ\'e \'egalement au rayon spectral $\rho$. 
\par
Si de plus $C$ est primitive, alors la suite de matrices $(\frac{1}{\rho^n}C^n)_{n\geq 1}$ converge vers la matrice
$G:=vw^T/w^Tv$ qui est la projection sur le noyau de $C-\rho I$ parall\`element \`a l'image de $C-\rho I$ (\cite{meyer}, p. 674). On observe que $G$  a tous ses coefficients strictement positifs.

\section{Preuves des th\'eor\`emes 2 et 3}

Consid\'erons d'abord un polyn\^ome $p(x)=x^k-c_{k-1}x^{k-1}-\cdots-c_1x-c_0$ comme dans la premi\`ere partie du th\'eor\`eme 1, c'est-\`a-dire tel que $c_j\geq 0$ pour tout $1\leq j\leq k-1$ et $c_0>0$. On lui associe sa \emph{matrice compagnon}
$$
C=\left(
\begin{array}{ccccc}
0 & 1 & 0 & \ldots & 0\\
0 & 0 & 1 & \ldots & 0\\
\vdots & \vdots & \ddots & \vdots & \vdots\\
0 & 0 & 0 & \ldots & 1\\
c_0 & c_1 & \ldots & c_{k-2} & c_{k-1}
\end{array}
\right).
$$
Cela signifie que le polyn\^ome caract\'eristique $p_C(x)=\det(C-xI)$ satisfait : $p_C(x)=(-1)^kp(x)$. Les deux polyn\^omes ont par cons\'equent exactement les m\^emes racines avec les m\^emes multiplicit\'es. Comme $c_0>0$, $C$ est irr\'eductible puisque, dans le graphe associ\'e, il y a une arête de 1 vers 2, de 2 vers 3, etc. de $k-1$ vers $k$ et de $k$ vers 1, au moins. Ainsi, $p(x)$ satisfait les conditions de la premi\`ere partie du th\'eor\`eme 2 si et seulement si sa matrice compagnon est primitive.
\par
La partie (1) du th\'eor\`eme 2 est alors une cons\'equence imm\'ediate de la proposition 7 ci-dessous, qui est un cas particulier du th\'eor\`eme de la page 679 de \cite{meyer}, mais nous en donnons une d\'emonstration par souci d'\^etre complet.

\begin{prop}
Soient $c_0,c_1,\ldots,c_{k-1}$ des nombres r\'eels positifs ou nuls, $c_0>0$ et soit $I=\{1\leq j<k:c_j>0\}$. Alors la matrice 
$$
C=\left(
\begin{array}{ccccc}
0 & 1 & 0 & \ldots & 0\\
0 & 0 & 1 & \ldots & 0\\
\vdots & \vdots & \ddots & \vdots & \vdots\\
0 & 0 & 0 & \ldots & 1\\
c_0 & c_1 & \ldots & c_{k-2} & c_{k-1}
\end{array}
\right).
$$
est primitive si et seulement si $\mathrm{pgcd}(I\cup\{k\})=1$.
\end{prop} 
\emph{Preuve.} Comme nous l'avons observ\'e ci-dessus, $C$ est irr\'eductible. Notons $N$ le cardinal de $I$.
Soit $S$ la matrice
$$
S=\left(
\begin{array}{ccccc}
0 & 1 & 0 & \ldots & 0\\
0 & 0 & 1 & \ldots & 0\\
\vdots & \vdots & \ddots & \vdots & \vdots\\
0 & 0 & 0 & \ldots & 1\\
\frac{1}{N+1} & \ldots & \frac{1}{N+1} & \ldots & \frac{1}{N+1}
\end{array}
\right)
$$
où dans la derni\`ere ligne, pour tout $j\in I\cup\{0\}$, $c_j$ est remplac\'e par $\frac{1}{N+1}$. La matrice ainsi obtenue est donc stochastique. En particulier, $\rho(S)=1$, et $S$ est irr\'eductible. De plus, par les rappels du paragraphe 2, $C$ est primitive si et seulement si $S$ l'est. On remplace alors $C$ par $S$, qui est la matrice compagnon de 
$$
q(x)=x^k-\sum_{j\in I\cup\{0\}}\frac{x^j}{N+1}.
$$
Supposons d'abord que $d=\mathrm{pgcd}(I\cup\{k\})>1$. \'Ecrivons $j=a_jd$ et $k=ad$ avec $a_j,a\in\N^*$ pour tout $j\in I$, et consid\'erons le polyn\^ome
$$
r(x)=x^a-\sum_{j\in I\cup\{0\}}\frac{x^{a_j}}{N+1}.
$$
Alors $q(x)=r(x^d)$,
$r(1)=0$, et pour tout $\zeta\in\C$, $\zeta^d=1$, on a $q(\zeta)=r(1)=0$. Cela d\'emontre que $C$ n'est pas primitive.
\par
R\'eciproquement, supposons que $\mathrm{pgcd}(I\cup\{k\})=1$, et soit $\zeta$ une racine de $q(x)$ telle que $|\zeta|=1$. Par le paragraphe 2, on sait que $\zeta$ est alors une racine de l'unit\'e. \'Ecrivons donc 
$$
\zeta=e^{2i\pi m/d}
$$
avec $1\leq m\leq d$ et $\mathrm{pgcd}(m,d)=1$. L'\'egalit\'e $q(\zeta)=0$ donne
$$
(N+1)e^{2i\pi km/d}=\sum_{j\in I}e^{2i\pi jm/d}+1,
$$
et en prenant les modules
$$
N+1=\left|\sum_{j\in I}e^{2i\pi jm/d}+1\right|\leq \sum_{j\in I}|e^{2i\pi jm/d}|+1=N+1
$$
qui implique que $d$ divise $jm$ pour tout $j\in I$, donc que $e^{2i\pi jm/d}=1$, puis que $e^{2i\pi km/d}=1$, donc que $d$ divise \'egalement $km$. Ces conditions impliquent que $d$ divise $\mathrm{pgcd}(I\cup\{k\})$ ; comme
$\mathrm{pgcd}(m,d)=1$, on obtient que $\zeta=1$, et $S$ est primitive.
\hfill $\square$

\par\vspace{3mm}
Pour d\'emontrer la seconde affirmation du th\'eor\`eme 2, consid\'erons un polyn\^ome \`a coefficients r\'eels $p(x)=x^k-cx^{k-1}-\sum_{j=0}^{k-2}c_jx^j$ tel que
$c\geq 2$ et $c>\sum_{j=0}^{k-2}|c_j|+1$.
On observe d'abord que pour tout $z\in\C$ tel que $|z|\geq 1$, on a 
$$
\left|\sum_{j=0}^{k-2}c_jz^j\right|<|z|^{k-2}(c-1).
$$
Cela implique facilement que $p(c-1)<0$, que $p(c+1)>0$ et que $p(z)\not=0$ pour tout $|z|=1$. Ainsi, $p$ s'annule en un nombre r\'eel $1\leq c-1<\rho<c+1$. Enfin, soit $q(z)=z^k-cz^{k-1}$. On a $|q(z)|\geq c-1$ pour tout $|z|=1$, et, pour ces m\^emes valeurs de $z$, on a
$$
|p(z)-q(z)|=\left|\sum_{j=0}^{k-2}c_jz^j\right|<c-1\leq |q(z)|.
$$
Puisque $q$ admet $k-1$ zéros (compt\'es avec multiplicit\'es) dans le disque unité ouvert, il en est de m\^eme pour $p$ par le th\'eor\`eme de Rouch\'e.
\hfill $\square$

\par\vspace{3mm}
Consid\'erons ensuite une suite $(a_n)_{n\geq 1}$ comme dans la troisième partie du théorème 2 :
$a_n>0$ pour tout $n$, et elle satisfait la relation de r\'ecurrence
$$
a_{n+k}=c_{k-1}a_{n+k-1}+c_{k-2}a_{n+k-2}+\ldots+c_1a_{n+1}+c_0a_n\quad (n\geq 1)
$$
où les $c_i$ satisfont l'une des deux premières conditions du théorème 2. On note encore $\rho$ l'unique racine positive du polyn\^ome caract\'eristique $p(x)=x^k-c_{k-1}x^{k-1}-\ldots-c_0$ de la relation de récurrence.
\par
Pour d\'emontrer la troisième partie du th\'eor\`eme 2, il suffit de v\'erifier que 
$\lim_{n\to\infty}\frac{a_n}{\rho^n}$ existe et est positive (\cite{jol2}, th\'eor\`eme 2.4). 
Le cas où les coefficients $c_i$ satisfont la condition (1) est une conséquence de la proposition suivante.

\begin{prop}
Avec les hypoth\`eses et les notations ci-dessus, si les $c_i$ satisfont la condition (1) du théorème 2, on a :
$$
\lim_{n\to\infty}\frac{a_n}{\rho^n}
$$ existe et est positive.
\end{prop}
\emph{Preuve.} Notons encore $C$ la matrice compagnon du polyn\^ome $p(x)$ et posons pour
tout entier $n\geq 1$ 
$$
A_n=\left(
\begin{array}{c}
a_n\\
a_{n+1}\\
\vdots\\
a_{n+k-1}
\end{array}\right)
$$
de sorte que 
$$
CA_n=A_{n+1}
$$
pour tout $n>0$. Par suite, on a $A_n=C^{n-1}A_1$ pour tout $n$. Si $v$ et $w$ sont des vecteurs de Perron-Frobenius pour $C$ et $C^T$ respectivement, on a vu que la suite de matrices $(\frac{1}{\rho^n}C^n)$ converge vers la matrice $G=vw^T/w^Tv$ lorsque $n\to\infty$. Ainsi,
$$
\frac{1}{\rho^n}A_n=\frac{1}{\rho}\cdot\frac{1}{\rho^{n-1}}C^{n-1}A_1\to_{n\to\infty}\frac{1}{\rho}GA_1.
$$
En particulier, $\frac{a_n}{\rho^n}$, qui est la premi\`ere composante de $\frac{1}{\rho^n}A_n$, converge vers la premi\`ere composante de $\frac{1}{\rho}GA_1$, et comme $A_1>0$ et $G>0$, toutes les composantes de $GA_1$ sont positives.
\hfill $\square$

\par\vspace{3mm}
Pour terminer la preuve du théorème 2, considérons une suite $(a_n)_{n\geq 1}\subset\R^+\setminus\{0\}$ qui satisfait la relation de récurrence 
$$
a_{n+k}=c_{k-1}a_{n+k-1}+\ldots+c_1a_{n+1}+c_0a_n
$$
avec $c_{k-1}\geq 2$, $c_{k-1}>\sum_{j=0}^{k-2}|c_j|+1$, et telle que $a_n$ ne converge pas vers 0 lorsque $n\to\infty$. 
D'après la structure de telles suites (rappelée dans le paragraphe 1), si $\zeta_1=\rho>1$, $\zeta_2,\ldots,\zeta_m$ sont les racines du polynôme 
caractéristique $p(x)=x^k-c_{k-1}x^{k-1}-\ldots-c_1x-c_0$, de multiplicités respectives $\mu_1=1$, $\mu_2,\ldots,\mu_m$, il existe des coefficients $\alpha_{j,\ell}$ tels que 
$$
a_n=\alpha_{1,1}\rho^n+\sum_{j=2}^m\sum_{\ell=0}^{\mu_j-1}\alpha_{j,\ell}n^\ell \zeta_j^n\quad \forall n\geq 1.
$$
Puisque $\rho>1$ et $|\zeta_j|<1$ pour tout $j=2,\ldots,m$, $a_n$ est de la forme $a_n=\alpha_{1,1}\rho^n+\beta_n$ avec $\beta_n\to 0$ lorsque $n\to\infty$.
Les hypothèses $a_n>0$ pour tout $n$ et $a_n\not\to 0$ impliquent que $\alpha_{1,1}>0$, et on obtient immédiatement que 
$\lim_{n\to\infty}\frac{a_n}{\rho^n}=\alpha_{1,1}.$
\hfill $\square$

\par\vspace{3mm}
Nous passons enfin \`a la preuve du th\'eor\`eme 3.

\par\vspace{3mm}\noindent
\emph{Preuve du th\'eor\`eme 3.} En n'\'ecrivant que les coefficients non nuls dans l'expression du polyn\^ome $p(x)$, on a
$$
p(x)=x^k-c_{k-k_1}x^{k-k_1}-\ldots-c_{k-k_s}x^{k-k_s}-c_0
$$
avec $1\leq k_1<\ldots<k_s<k$. Par le th\'eor\`eme de la page 679 de \cite{meyer}, on obtient :
$$
h=\mathrm{pgcd}(k-k_1,\ldots,k-k_s,k)=\mathrm{pgcd}(k_1,\ldots,k_s,k)=\mathrm{pgcd}(I\cup\{k\}).
$$
\'Ecrivons, pour $1\leq i\leq s$, $k_i=hk'_i$, et aussi $k=hk'$. Pour chaque $m$ fix\'e, en rempla\c{c}ant $n$ par $m+hn$ dans la relation de r\'ecurrence, on obtient
$$
a_{m+h(n+k')}=c_{k-k_1}a_{m+h(n+k'-k'_1)}+\ldots+c_{k-k_s}a_{m+h(n+k'-k'_s)}+c_0a_{m+hn}.
$$
Cela signifie que la suite $(a_{m+hn})_{n\geq 1}$ satisfait une relation de r\'ecurrence d'ordre $k'$ dont le polyn\^ome caract\'eristique $q(x)$ est 
$$
q(x)=x^{k'}-c_{k-k_1}x^{k'-k'_1}-\ldots-c_{k-k_s}x^{k'-k'_s}
$$
avec $\mathrm{pgcd}(I(q)\cup\{k'\})=1$. On applique alors les conclusions des propositions pr\'ec\'edentes \`a la suite $(a_{m+hn})_{n\geq 1}$.
\hfill $\square$

\section{Crit\`ere d'irrationnalit\'e de $\log_b(\rho)$}

Soit $p(x)$ le polyn\^ome $p(x)=x^k-c_{k-1}x^{k-1}-\ldots-c_1x-c_0$ avec $c_i\in\R^+$ et $c_0>0$. On d\'esigne encore par $I$ l'ensemble des indices $1\leq i<k$ tels que $c_i>0$,
et par $\rho$ la racine positive de $p(x)$. On suppose encore que $I\not=\emptyset$.
\par\vspace{3mm}
Le r\'esultat suivant sera utilis\'e dans la preuve du th\'eor\`eme 6.

\begin{prop}
Soit $b>2$ un entier. On suppose que les $c_i\in\Q^+$ et que :
\begin{enumerate}
	\item [(1)] $\mathrm{pgcd}(I\cup\{k\})=1$ ;
	\item [(2)] $\rho$ n'est ni entier, ni inverse d'entier.
\end{enumerate}
Alors $\log_b(\rho)$ est irrationnel.
\end{prop}
\emph{Preuve.} Supposons que $\log_b(\rho)$ soit rationnel, de sorte qu'il existe deux entiers $p$ et $q>0$ premiers entre eux tels que $\rho=b^{p/q}$. Comme $\rho\not=1$, on a $p\not=0$. D\'esignons par $m(x)$ le polyn\^ome minimal (unitaire) de $\rho$.
\par
Si $\rho$ est irrationnel, alors $m(x)$ est de degr\'e au moins 2 et $\rho$ en est une racine simple. 
Mais $\rho$ est également une racine de $t(x)=x^q-b^p$, donc le polyn\^ome $m(x)$ divise $t(x)$, et ainsi toutes les racines de $m(x)$ sont des racines de $t(x)$. Or, ces derni\`eres sont toutes de module \'egal \`a $\rho$. Cela contredit (1) car $m(x)$ divise \'egalement $p(x)$, et ce dernier aurait plusieurs racines de module \'egal \`a $\rho$. 
\par
Par suite, $\rho$ est n\'ecessairement rationnel ; il existe des entiers positifs $\alpha,\beta$ tels que $\rho=\frac{\alpha}{\beta}$ et $\mathrm{pgcd}(\alpha,\beta)=1$. On obtient $\alpha^q=b^p\beta^q$, et l'unicit\'e de la d\'ecomposition en facteurs premiers implique que
\begin{itemize}
	\item [] $\beta=1$ si $p$ est positif, c'est-à-dire si $\rho>1$, et alors $\rho=\alpha$ est entier, 
	\item [] $\alpha=1$ si $p$ est négatif, c'est-à-dire si $\rho<1$, et alors $\rho=1/\beta$ est un inverse d'entier, 
\end{itemize}
ce qui contredit la seconde hypoth\`ese.
\hfill $\square$
\par\vspace{3mm}\noindent
\textbf{Remarque.} La condition (2) est facile \`a utiliser puisqu'il est commode de localiser $\rho$ dans $\R^+$ gr\^ace à une étude succinte du polyn\^ome $p(x)$.

\begin{exemp}
Soit $m\geq 1$ un entier fixé ; pour tout entier $k\geq 2$, soit 
$$
p_{k,m}(x)=x^k-mx^{k-1}-mx^{k-2}-\ldots-mx-m,
$$
qui généralise le polyn\^ome caractéristique des suites de Fibonacci d'ordre $k$ introduites dans l'exemple 4, et qui satisfait les conditions du th\'eor\`eme 2. Notons $\rho_{k,m}$ la racine positive de $p_{k,m}(x)$. 
\par
Elle est irrationnelle par la proposition 5, $\rho_{k,m}>1$ car $p_{k,m}(1)=1-km<0$, et en fait, $\rho_{k,m}$ est un \emph{nombre de Pisot} (cf \cite{brauer}, \cite{pisot2}) : c'est un entier alg\'ebrique, et $p_{k,m}(x)$ est son polyn\^ome minimal car on a 
$c_{k-1}\geq\ldots\geq c_0>0$, et toutes les racines $\zeta\not=\rho_{k,m}$ de $p_{k,m}(x)$ satisfont $|\zeta|<\rho_{k,m}$.
\par
On a $m<\rho_{k,m}<m+1$ car 
$$
p_{k,m}(m)=
\left\{
\begin{array}{ll}
1-k<0, & m=1,\\
\frac{m-m^k}{m-1}<0, & m\geq 2,
\end{array}\right.
$$
et $p_{k,m}(m+1)=1$ pour tout $m$.
\par
Plus pr\'ecis\'ement, on va d\'emontrer que $\rho_{k,m}<\rho_{k+1,m}$ pour tout $k$ assez grand et pour tout $m$, et que
$$
\frac{(m+1)k}{k+1}<\rho_{k,m}<m+1.
$$
Cela d\'emontrera que, pour tout $m$ fix\'e, la suite $(\rho_{k,m})$ est croissante et converge vers $m+1$.
\par
La premi\`ere affirmation provient des égalités : 
$$
p_{k+1,m}(\rho_{k,m})=\rho_{k,m}^{k+1}-m\rho_{k,m}^k-\ldots-m\rho_{k,m}-m=\rho_{k,m}\cdot p_{k,m}(\rho_{k,m})-m=-m
$$
et $p_{k+1,m}(\rho_{k+1,m})=0$.
\par
Pour d\'emontrer la seconde affirmation, nous introduisons le polyn\^ome auxiliaire 
$$
q(x)=(x-1)p_{k,m}(x)=x^{k+1}-(m+1)x^k+m.
$$
On v\'erifie sans peine que
$$
q\left(\frac{(m+1)k}{k+1}\right)=-\frac{m+1}{k+1}\left(\frac{(m+1)k}{k+1}\right)^k+m.
$$
Or, l'inégalité $-\frac{m+1}{k+1}\left(\frac{(m+1)k}{k+1}\right)^k+m<0$ est équivalente à $(m+1)^{k+1}>m(k+1)(1+1/k)^k$ qui est vraie d\`es que $k$ est assez grand. Cela d\'emontre la seconde affirmation.
\end{exemp}

\par\vspace{3mm}\noindent
\textbf{Remarques.} (1) Les suites $(\rho_{k,m}^\ell)_{\ell\geq 1}$ ne sont pas \'equidistribu\'ees $\bmod\ 1$ car la distance entre $\{\rho_{k,m}^\ell: \ell\leq n\}$ et $\N$ tend vers 0. Mais qu'en est-il de $(\alpha\cdot\rho_{k,m}^\ell)_{\ell\geq 1}$ pour 
$\alpha$ irrationnel positif ? (On sait que, pour presque tout nombre réel $x>1$, la suite $(x^n)_{n\geq 1}$ est équidistribuée
$\bmod\ 1$ ; cf \cite{KN}, chap.1, corollaire 4.2.)\\
(2) Si $c_0,c_1,\ldots,c_{k-1}\in\Z$ avec $k\geq 2$ sont tels que 
$$
c_{k-1}>\sum_{j=0}^{k-2}|c_j|+1
$$
et $c_0\not=0$, le polyn\^ome associ\'e $p(x)=x^k-c_{k-1}x^{k-1}-\ldots-c_0$ est le polyn\^ome minimal d'un nombre de Pisot $\rho_k$ (\cite{pisot2}, chap. 5.2), et cela donne une famille de suites qui satisfont une relation de r\'ecurrence \`a coefficients entiers non n\'ecessairement positifs et qui satisfont \'egalement
la loi de Benford d'après le théorème 2.
\par\vspace{3mm}
Nous passons enfin \`a la preuve du th\'eor\`eme 6 :

\par\vspace{3mm}\noindent
\emph{Preuve du th\'eor\`eme 6.} (1) Pour $N\geq 3$ fix\'e, posons
$$
B_N=\{2<b\leq N: \log_b(\rho)\in \Q\}.
$$
Il suffit de d\'emontrer que 
$$
|B_N|\leq\frac{\sqrt N\log(N)}{\log(2)}.
$$
C'est \'evident si $B_N$ est vide ou s'il ne contient qu'un \'el\'ement. S'il contient au moins deux \'el\'ements, soit $b_0=\min\{b\in B_N\}$. Pour tout $b\in B_N$ tel que $b>b_0$, en utilisant encore la décomposition en facteurs premiers, on vérifie qu'
il existe un entier $u\leq\sqrt N$ tel que $b_0$ et $b$ soient des puissances enti\`eres de $u$. Par suite,
$$
B_N\subset\{2<u^p\leq N:p\geq 1, u\leq\sqrt N\}.
$$
On obtient ainsi la majoration annonc\'ee.
\par\noindent
(2) Soit $h\geq 1$ le nombre de racines positives distinctes du polyn\^ome $p(x)$. 
Les sous-suites qui ont des conditions initiales positives parmi $(a_{m+hn})_{n\geq 1}$, $0\leq m<h$, ont un polyn\^ome caract\'eristique dont $\rho^{1/h}$ est la racine positive. Si $\rho$ n'est ni entier ni inverse d'entier, alors $\rho^{1/h}$ non plus et les conditions de la proposition 9 sont satisfaites.
\par
Enfin, si $\rho$ est rationnel, en utilisant la d\'ecomposition en facteurs premiers, on v\'erifie comme ci-dessus que $b^{\pm p/q}\in\N$. Par suite,
il existe des entiers $c,d\geq 1$ et $u\geq 2$ tels que $\rho^{\pm 1}=u^c$ et $b=u^d$, mais alors $\rho$ serait entier ou inverse d'un entier, contrairement \`a l'hypoth\`ese.
Cela d\'emontre la deuxi\`eme partie du th\'eor\`eme 6.
\hfill $\square$

\par\vspace{3mm}

\bibliographystyle{plain}
\bibliography{refprim}
\par
\vspace{1cm}\noindent
H. Deligny\\
Professeur agr\'eg\'e de Math\'ematiques,
Acad\'emie de La R\'eunion,
20 rue Colbert,
Saint Paul La R\'eunion,
\texttt{abozinis@hotmail.com}

\vspace{3mm}\noindent
P. Jolissaint\\
 Institut de Math\'emathiques,
       Universit\'e de Neuch\^atel,
       Emile-Argand 11,
       CH-2000 Neuch\^atel,
       \texttt{paul.jolissaint@unine.ch}

\end{document}